\documentclass[12pt, reqno]{amsart}
\usepackage{amsmath, amsthm, amscd, amsfonts, amssymb, graphicx, color}
\usepackage[bookmarksnumbered, colorlinks, plainpages]{hyperref}
\hypersetup{colorlinks=true,linkcolor=red, anchorcolor=green, citecolor=cyan, urlcolor=red, filecolor=magenta, pdftoolbar=true}

\newtheorem{theorem}{Theorem}[section]
\newtheorem{lemma}[theorem]{Lemma}

\theoremstyle{definition}

\theoremstyle{remark}

\numberwithin{equation}{section}

\begin{document}

\setcounter{page}{1}

\title[Non-squareness and local uniform non-squareness]{Non-squareness and local uniform non-squareness properties of Orlicz-Lorentz function spaces endowed with the Orlicz norm}

\author[B. Chen \MakeLowercase{and} W. Gong]{Bowen Chen \MakeLowercase{and} Wanzhong Gong$^{*}$}

\address{Department of Mathematics, Anhui Normal University, Wuhu 241000, China.}
\email{\textcolor[rgb]{0.00,0.00,0.84}{675024965@qq.com; gwzsjy@ahnu.edu.cn}}




\keywords{non-squareness; Orlicz-Lorentz function space; Orlicz norm.\newline
\indent $^{*}$Corresponding author}


\begin{abstract}
In this paper the necessary and sufficient conditions were given for Orlicz-Lorentz function space endowed with the Orlicz norm having non-squareness and local uniform non-squareness. 
\end{abstract} 

\maketitle

\section{Introduction}  
Let $X$ be a Banach space, We say that a Banach space $X$ is uniformly non-square if there exists $\delta\in(0,1)$ such that
 \[\min\left\{\left\|\frac{x-y}{2}\right\|,\left\|\frac{x+y}{2}\right\|\right\}\leq1-\delta\]
 for any $x,y\in S(X)$$^\text{\cite{James-1}}$. We have known that uniform non-squareness implies super-reflexivity$^\text{\cite{James-2}}$ and the fixed point property$^\text{\cite{Garcia}}$. In 1985 Hudzik gave the criterion for Orlicz space being  uniformly non-square$^\text{\cite{Hudzik-2}}$, then Hudzik and Cui etc obtained the characteristic for Orlicz-Bochner space being uniformly non-square$^\text{\cite{Hudzik-1,Shang}}$.
 We say that a Banach space $X$ is non-square if
 \[\min\left\{\left\|\frac{x-y}{2}\right\|,\left\|\frac{x+y}{2}\right\|\right\}<1\]
 for any $x,y\in S(X)$$^\text{\cite{James-1}}$. A Banach space $X$ is said to be locally uniformly non-square if for any $x\in S(X)$ there exists $\delta=\delta(x)\in(0,1)$ such that
 \[\min\left\{\left\|\frac{x-y}{2}\right\|,\left\|\frac{x+y}{2}\right\|\right\}\leq1-\delta\]
 for any $y\in B(X)$. In 2013 Foralewski, Hudzik and Kolwicz show the sufficient and necessary conditions for Orlicz-Lorentz space with the Luxemburg norm being non-square and uniformly non-square$^\text{\cite{Foralewski,Foralewski-1}}$. Recently Foralewski and Ko\'{n}czak got the  criterion for Orlicz-Lorentz space with the Luxemburg norm being locally uniformly non-square$^\text{\cite{Foralewski-2}}$.
   In this paper we will discuss the characteristic of Orlicz-Lorentz function space with the Orlicz norm being non-square and locally uniformly non-square. For more reference about non-squareness and local uniform non-squareness we refer to \cite{Hudzik-1,Hudzik-2,Hudzik-3,Kondagunta,Prus}.

A function $\varphi$ : $\mathbb{R}\rightarrow\mathbb{R_+}$ is said to be an Orlicz function$^\text{\cite{Krasnoselskii}}$  if $\varphi$ is convex, even, $\varphi(0)=0$, $\varphi(u)>0$ for all $u>0$, $\underset{u\rightarrow0}{\lim}\frac{\varphi(u)}{u}=0$ and $\underset{u\rightarrow\infty}{\lim}\frac{\varphi(u)}{u}=\infty$. Its {\it complementary function} $\psi$ is defined by
\begin{align*}
\psi(v)=\underset{u>0}{\sup}\{|uv|-\varphi(u)\}
\end{align*}
for all $v\in\mathbb{R}$.
If $\varphi$ is an Orlicz function, then its complementary function $\psi$ is an Orlicz function. Recall that an Orlicz function $\varphi$ satisfies the condition $\Delta_2$ for all value ($\varphi\in\Delta_2(\mathbb{R})$ for short) if there exists a constant $K>0$ such that $\varphi(2u)\leq K\varphi(u)$ holds for all $u\in\mathbb{R}$. Analogously, an Orlicz function $\varphi$ satisfies the condition $\Delta_2$ for large values ($\varphi\in\Delta_2(\infty)$ for short) if there exist a constant $K>0$ and a constant $u_0\geq0$ such that $\varphi(2u)\leq K\varphi(u)$ holds for all $u\geq u_0$. Sometimes we say $\varphi$ satisfy $\nabla_2$ condition for all values (or $\nabla_2$ condition for large values) if  $\psi$ satisfy $\Delta_2$ condition for all values (or $\Delta_2$ condition for large values).  For more properties of Orlicz function, we may refer to \cite{chen-2,Krasnoselskii}.

For any measurable function $x$ : $[0,\gamma)\rightarrow\mathbb{R}$,  where $\gamma\leq\infty$,  its {\it distribution function}  and {\it decreasing rearrangement}  are defined as follows
\begin{align*}
&d_x(\theta)=\mu\{s\in[0,\gamma):|x(s)|>\theta\},\\
&x^\ast(t)=\inf\{\theta>0:d_x(\theta)\leq t\},\quad t\geq0,
\end{align*}
where $\mu$ denotes the Lebesgue measure. A function $\omega$ : $[0,\gamma)\rightarrow\mathbb{R_+}$ is said to be a {\it weight function} if it is non-increasing and locally integrable. Define $\alpha=\sup\{t\geq0: \omega(t)>0\}$, $W(x)=\int_{0}^{x}\omega(s)\,ds$.

The {\it Orlicz-Lorentz space} $\Lambda_{\varphi,\omega}$ is the set of all Lebesgue measurable functions $x$ on $[0,\gamma)$ such that
\begin{align*}
\rho_{\varphi,\omega}(\lambda x)=\int_{0}^{\gamma} \varphi(\lambda x^\ast(t))\omega(t)\,dt<\infty
\end{align*}
for some $\lambda>0$. It is known that the Orlicz-Lorentz space endowed with the Luxemburg norm
\begin{align*}
\|x\|_{\varphi,\omega}=\inf\left\{\varepsilon>0:\rho_{\varphi,\omega}\Big(\frac{x}{\varepsilon}\Big)\leq1\right\}
\end{align*}
is a Banach space$^\text{\cite{Kaminska-1}}$.
If $\varphi(t)=t$, then $\Lambda_{\varphi,\omega}$ is the Lorentz function space $L_{1,\omega}$. The $L_{1,\omega}${\it -norm} of $x\in L_{1,\omega}$ is defined by
\begin{align*}
\|x\|_{1,\omega}=\int_{0}^{\gamma} x^\ast(t)\omega(t)\,dt.
\end{align*}
Obviously,
\begin{align*}
\rho_{\varphi,\omega}(x)=\|\varphi\circ x\|_{1,\omega}.
\end{align*}

Recall a Banach lattice $E=(E,\leq,\|\cdot\|)$ is said to be strictly monotone$^\text{\cite{Birkhoff}}$ if $x,y\in E$, $0\leq y\leq x$ and $y\neq x$ imply that $\|y\|<\|x\|$. A Banach  lattice $E=(E,\leq,\|\cdot\|)$ is said to be lower locally uniformly monotone$^\text{\cite{Hudzik-k-1}}$, whenever for any $x\in(E)_{+}$ (the positive cone of $E$) with $\|x\|=1$ and any $\varepsilon\in(0,1)$ there exists $\delta(x,\varepsilon)\in(0,1)$ such that the conditions $0\leq y\leq x$ and $\|y\|\geq\varepsilon$ imply that $\|x-y\|\leq1-\delta(x,\varepsilon)$.
 Suppose $A_1,A_2$ are the subsets of $\mathbb{R}$. A mapping $\sigma$ : $A_1\rightarrow A_2$ is called {\it measure preserving transformation} if for any measurable set $E\subset A_2$, it holds that the set $\sigma^{-1}(E)\subset A_1$ is also measurable and $\mu(\sigma^{-1}(E))=\mu(E)$.
From \cite{Bennett}  we know that if $x$ is a simple function with compact support in $\Lambda_{\varphi,\omega}[0,\gamma)$ then there is a measure preserving transformation $\sigma$ such that
\begin{align*}
\int_{0}^{\gamma}\varphi(x^\ast)\omega=\int_{0}^{\gamma}\varphi(x)\omega\circ\sigma.
\end{align*}

In 1999, Wu and Ren defined the Orlicz norm on the space $\Lambda_{\varphi,\omega}[0,\gamma)$ for $\gamma<\infty$$^\text{\cite{Wu}}$ as
\begin{align*}
\|x\|^\circ_{\varphi,\omega}=\underset{\rho_{\psi,\omega}(y)\leq1}{\sup}\int_{0}^{\gamma}x^\ast(t)y^\ast(t)\omega(t)\,dt.
\end{align*}
In \cite{Wu} the authors proved  that endowed with the Orlicz norm, $\Lambda_{\varphi,\omega}[0,\gamma)$ is a Banach space (denoted by $\Lambda^\circ_{\varphi,\omega}[0,\gamma)$), and obtained the following properties of $\Lambda^\circ_{\varphi,\omega}[0,\gamma)$:

(i)\ Let $x\in\Lambda_{\varphi,\omega}[0,\gamma)$. If $\|x\|^\circ_{\varphi,\omega}\leq1$, then $\rho_{\varphi,\omega}(x)\leq\|x\|^\circ_{\varphi,\omega}$.

(ii)\ For any $x\in\Lambda_{\varphi,\omega}[0,\gamma)$, $\|x\|_{\varphi,\omega}\leq\|x\|^\circ_{\varphi,\omega}\leq2\|x\|_{\varphi,\omega}$.

(iii)\ If there exists $k>1$ such that $\rho_{\psi,\omega}(p(k|x|))=1$, then
\begin{align*}
\|x\|^\circ_{\varphi,\omega}=\int_{0}^{\gamma} x^\ast(t)p(kx^\ast(t))\omega(t)\,dt=\frac{1}{k}(1+\rho_{\varphi,\omega}(kx)).
\end{align*}

(iv)\ $\|x\|^\circ_{\varphi,\omega}=\underset{\rho_{\psi,\omega}(y)\leq1}{\sup}\int_{0}^{\gamma}x^\ast(t)y^\ast(t)\omega(t)\,dt=\underset{k>0}{\inf}\frac{1}{k}(1+\rho_{\varphi,\omega}(kx))$.

(v)\ For any $x\in\Lambda_{\varphi,\omega}[0,\gamma)$,
\begin{align*}
\|x\|^\circ_{\varphi,\omega}=\frac{1}{k_x}(1+\rho_{\varphi,\omega}(k_xx))
\end{align*}
if and only if $k_x\in K(x)=[k^\ast,k^{\ast\ast}]$,
where $k^\ast=\inf\{h>0:\rho_{\psi,\omega}(p(h|x|))\geq1\}$ and $k^{\ast\ast}=\sup\{h>0:\rho_{\psi,\omega}(p(h|x|))\leq1\}$. For the sake of convenience, in this paper we will consider $\gamma=1$ whenever $\gamma<+\infty$.

\section{Some Lemmas}
In 2012, Wang and Chen extended the definition and properties of \cite{Wu} to $\gamma\leq\infty$$^\text{\cite{Wang-2}}$,
and get the following two lemmas.

\begin{lemma}\label{1-le}
{\rm(I)}\ $\inf\{k:k\in K(x), \|x\|^\circ_{\varphi,\omega}=1\}>1$ if and only if $\varphi\in\Delta_2$.

{\rm(II)}\ The set $Q=\bigcup\{K(x): a\leq\|x\|^\circ_{\varphi,\omega}\leq b\}$ is bounded for each $b\geq a>0$ if and only if $\varphi\in\nabla_2$.
\end{lemma}

\begin{lemma}[\cite{Wang-2}]\label{2-le}
Let $A\subset[0,\infty)$ and $\mu A=t$.

{\rm(I)} For $t<\infty$,
\begin{align*}
\|\chi_A\|^\circ_{\varphi,\omega}=\psi^{-1}\left(\frac{1}{W(t)}\right)W(t);
\end{align*}

{\rm(II)} For $t=\infty$ and $W(\infty)<\infty$,
\begin{align*}
\|\chi_A\|^\circ_{\varphi,\omega}=\psi^{-1}\left(\frac{1}{W(\infty)}\right)W(\infty).
\end{align*}
\end{lemma}

By the same method as the proof of Theorem 2 in \cite{Kaminska-2}, we can get that if $\psi\notin\Delta_2(\infty)\ ({\rm or}\ \psi\notin\Delta_2(\mathbb{R}))$, then for any $\varepsilon\in(0,1)$, there exists a sequence $\{e_n\}$ with $\rho_{\varphi,\omega}(e_n)=\varepsilon$  such that
\begin{align*}
\frac{1}{8}\varepsilon\leq\left\|\sum_{i=1}^{n}a_ie_i\right\|^\circ_{\varphi,\omega}\leq1+\varepsilon,
\end{align*}
where $a_i\geq0$ and $\sum_{i=1}^{n}a_i=1$. Therefore we have
\begin{lemma}\label{5-le}
If $\psi\notin\Delta_2(\infty)\ ({\rm or}\ \psi\notin\Delta_2(\mathbb{R}))$, then $ \Lambda^\circ_{\varphi,\omega}[0,1)$ $({\rm or}\ \Lambda^\circ_{\varphi,\omega}[0,\infty))$ contains $\ell_1$.
\end{lemma}
%
%
%
%
%

\begin{lemma}[\cite{Hudzik,Kolwicz}]\label{6-le}
The Lorentz function space $L_{1,\omega}$ is strictly monotone if and only if $\omega$ is positive on $[0,\gamma)$ and $\int_{0}^{\gamma}\omega(t)\,dt=\infty$ whenever $\gamma=\infty$.
\end{lemma}

\begin{lemma}[\cite{Foralewski-2}]\label{7-le}
The Lorentz function space $L_{1,\omega}$ is lower locally uniformly monotone if and only if $\omega$ is positive on $[0,\gamma)$ and $\int_{0}^{\gamma}\omega(t)\,dt=\infty$ whenever $\gamma=\infty$.
\end{lemma}

\section{Non-squareness of $\Lambda^\circ_{\varphi,\omega}$}   
\begin{theorem}\label{3-th}
Orlicz-Lorentz space $\Lambda^\circ_{\varphi,\omega}[0,\infty)$ is non-square if and only if $\int_{0}^{\infty}\omega(t)\,dt=\infty$.
\end{theorem}
\proof (Necessity) If $\int_{0}^{\infty}\omega(t)\,dt<\infty$, let $A\subset[0,\infty)$ and $\mu A=\infty$. By Lemma \ref{2-le}, we obtain
\begin{align*}
\|\chi_A\|^\circ_{\varphi,\omega}=\psi^{-1}\left(\frac{1}{W(\infty)}\right)W(\infty).
\end{align*}
Let
\begin{align*}
x=&\frac{1}{\psi^{-1}\big(\frac{1}{W(\infty)}\big)W(\infty)}\chi_{A},\\
y=&\frac{1}{\psi^{-1}\big(\frac{1}{W(\infty)}\big)W(\infty)}\chi_{A_1}-\frac{1}{\psi^{-1}\big(\frac{1}{W(\infty)}\big)W(\infty)}\chi_{A_2},
\end{align*}
where $A_1\cup A_2= A$, $A_1\cap A_2=\emptyset$ and $\mu(A_1)=\mu(A_2)=\infty$. Since $x^\ast=y^\ast$, we get
\begin{align*}
\|x\|^\circ_{\varphi,\omega}=\|y\|^\circ_{\varphi,\omega}=\left\|\frac{x+y}{2}\right\|^\circ_{\varphi,\omega}=\left\|\frac{x-y}{2}\right\|^\circ_{\varphi,\omega}=1.
\end{align*}
Which is a contradiction.

(Sufficiency) Let $x,y\in S(\Lambda^\circ_{\varphi,\omega}[0,\infty))$. Fix $k_1\in K(x)$, $k_2\in K(y)$. Let $k=\frac{2k_1k_2}{k_1+k_2}$.

Denote
\begin{align*}
&A_1=\{t\in[0,\infty): x(t)y(t)>0\},\\
&A_2=\{t\in[0,\infty): x(t)y(t)<0\},\\
&A_3=\{t\in[0,\infty): x(t)y(t)=0 \ and \ \max\{|x(t)|,|y(t)|\}>0\}.
\end{align*}
If $t\in A_1$, we can get
\begin{align*}
\varphi\left(\frac{k(x-y)}{2}(t)\right)&=\varphi\left(\frac{k_2}{k_1+k_2}k_1x(t)-\frac{k_1}{k_1+k_2}k_2y(t)\right)\\
&<\varphi\left(\frac{k_2}{k_1+k_2}k_1x(t)+\frac{k_1}{k_1+k_2}k_2y(t)\right)\\
&\leq\frac{k_2}{k_1+k_2}\varphi(k_1x(t))+\frac{k_1}{k_1+k_2}\varphi(k_2y(t)).
\end{align*}
If $t\in A_2$, there holds
\begin{align*}
\varphi\left(\frac{k(x+y)}{2}(t)\right)&=\varphi\left(\frac{k_2}{k_1+k_2}k_1x(t)+\frac{k_1}{k_1+k_2}k_2y(t)\right)\\
&<\varphi\left(\frac{k_2}{k_1+k_2}k_1x(t)-\frac{k_1}{k_1+k_2}k_2y(t)\right)\\
&\leq\frac{k_2}{k_1+k_2}\varphi(k_1x(t))+\frac{k_1}{k_1+k_2}\varphi(k_2y(t)).
\end{align*}
If $t\in A_3$, the inequality $\varphi(\beta x)<\beta\varphi(x)$ for any $0<\beta<1$ follows that
\begin{align*}
\varphi\left(\frac{k(x+y)}{2}(t)\right)&=\varphi\left(\frac{k_2}{k_1+k_2}k_1x(t)+\frac{k_1}{k_1+k_2}k_2y(t)\right)\\
&<\max\left\{\frac{k_2}{k_1+k_2}\varphi(k_1x(t)),\frac{k_1}{k_1+k_2}\varphi(k_2y(t))\right\}\\
&\leq\frac{k_2}{k_1+k_2}\varphi(k_1x(t))+\frac{k_1}{k_1+k_2}\varphi(k_2y(t)).
\end{align*}
By the strict monotonicity of $L_{1,\omega}[0,\infty)$ we get, if $\mu(A_1)>0$,
\begin{align*}
\left\|\varphi\left(\frac{k(x-y)}{2}\right)\right\|_{1,\omega}&<\frac{k_2}{k_1+k_2}\|\varphi(k_1x)\|_{1,\omega}+\frac{k_1}{k_1+k_2}\|\varphi(k_2y)\|_{1,\omega}\\
&\leq\frac{k_2}{k_1+k_2}(k_1\|x\|^\circ_{\varphi,\omega}-1)+\frac{k_1}{k_1+k_2}(k_2\|y\|^\circ_{\varphi,\omega}-1)\\
&\leq k-1.
\end{align*}
Whence
\begin{align*}
\left\|\frac{x-y}{2}\right\|^\circ_{\varphi,\omega}\leq\frac{1+\rho_{\varphi,\omega}(\frac{k(x-y)}{2})}{k}<1.
\end{align*}
Similarly, if $\mu(A_2\cup A_3)>0$, there holds
\begin{align*}
\left\|\frac{x+y}{2}\right\|^\circ_{\varphi,\omega}<1.
\end{align*}
In summary
\begin{align*}
\min\left\{\left\|\frac{x+y}{2}\right\|^\circ_{\varphi,\omega}, \left\|\frac{x-y}{2}\right\|^\circ_{\varphi,\omega}\right\}<1.
\end{align*}

\begin{theorem}\label{4-th}
Orlicz-Lorentz space $\Lambda^\circ_{\varphi,\omega}[0,1)$ is non-square if and only if $\alpha:=\sup\{t\geq0:\omega(t)>0\}\in(\frac{1}{2},1]$.
\end{theorem}
\proof(Necessity) If $\alpha\in(0,\frac{1}{2}]$,
set
\begin{align*}
x=&\frac{1}{\psi^{-1}\big(\frac{1}{W(\alpha)}\big)W(2\alpha)}\chi_{[0,2\alpha)},\\
y=&\frac{1}{\psi^{-1}\big(\frac{1}{W(\alpha)}\big)W(2\alpha)}\chi_{[0,\alpha)}-\frac{1}{\psi^{-1}\big(\frac{1}{W(\alpha)}\big)W(2\alpha)}\chi_{[\alpha,2\alpha]}.
\end{align*}
Clearly,
\begin{align*}
\|x\|^\circ_{\varphi,\omega}=\|y\|^\circ_{\varphi,\omega}=\left\|\frac{x+y}{2}\right\|^\circ_{\varphi,\omega}=\left\|\frac{x-y}{2}\right\|^\circ_{\varphi,\omega}=1.
\end{align*}
Which shows that $\Lambda^\circ_{\varphi,\omega}[0,1)$ is not non-square, a contradiction.

(Sufficiency) Let $x,y\in S(\Lambda^\circ_{\varphi,\omega}[0,1))$. Fix $k_1\in K(x)$, $k_2\in K(y)$. Let $k=\frac{2k_1k_2}{k_1+k_2}$.

Case 1. Let $\alpha=1$. It is similar to the proof of Theorem \ref{3-th}.

Case 2. Suppose that $\frac{1}{2}<\alpha<1$. Define
\begin{align*}
A_{x,y}=\{t\in[0,1): \max\{|x(t)|,|y(t)|\}>0\}.
\end{align*}

Case 2.1 If $\mu(A_{x,y})\leq\alpha$. We define
\begin{align*}
\overset{\sim}x=x\chi_{A_{x,y}}\circ\sigma,\\
\overset{\sim}y=y\chi_{A_{x,y}}\circ\sigma,
\end{align*}
where $\sigma:[0,\mu(A_{x,y}))\rightarrow A_{x,y}$ is a measure preserving transformation. Obviously, $\varphi(\overset{\sim}x)$, $\varphi(\overset{\sim}y)$, $\varphi\left(\frac{\overset{\sim}x+\overset{\sim}y}{2}\right)$ and $\varphi\left(\frac{\overset{\sim}x-\overset{\sim}y}{2}\right)$ are equimeasurable with $\varphi(x\chi_{A_{x,y}})$, $\varphi(y\chi_{A_{x,y}})$, $\varphi\left(\frac{x+y}{2}\chi_{A_{x,y}}\right)$ and $\varphi\left(\frac{x-y}{2}\chi_{A_{x,y}}\right)$. Similarly as Case 1, we see
\begin{align*}
&\min\left\{\rho_{\varphi,\omega}\left(\frac{k(x-y)}{2}\right), \rho_{\varphi,\omega}\left(\frac{k(x+y)}{2}\right)\right\}\\
=&\min\left\{\rho_{\varphi,\omega}\left(\frac{k(x-y)}{2}\chi_{A_{x,y}}\right), \rho_{\varphi,\omega}\left(\frac{k(x+y)}{2}\chi_{A_{x,y}}\right)\right\}\\
=&\min\left\{\rho_{\varphi,\omega}\left(\frac{k\big(\overset{\sim}x-\overset{\sim}y\big)}{2}\right), \rho_{\varphi,\omega}\left(\frac{k\big(\overset{\sim}x+\overset{\sim}y\big)}{2}\right)\right\}\\
<&k-1.
\end{align*}
Therefore
\begin{align*}
\min\left\{\left\|\frac{x+y}{2}\right\|^\circ_{\varphi,\omega}, \left\|\frac{x-y}{2}\right\|^\circ_{\varphi,\omega}\right\}<1.
\end{align*}

Case 2.2 Assume $\mu(A_{x,y})>\alpha$. By the convexity of $\varphi$, we obtain
\begin{align*}
\varphi\left(\left(\frac{k(x+y)}{2}\right)^\ast(t)\right)&=\left(\varphi\left(\frac{k(x+y)}{2}\right)\right)^\ast(t)\notag\\
&=\left(\varphi\left(\frac{k_2}{k_1+k_2}k_1x+\frac{k_1}{k_1+k_2}k_2y\right)\right)^\ast(t)\notag\\
&\leq\left(\frac{k_2}{k_1+k_2}\varphi(k_1x)+\frac{k_1}{k_1+k_2}\varphi(k_2y)\right)^\ast(t).
\end{align*}
Similarly, for any $t\in[0,1)$,
\begin{align*}
\varphi\left(\left(\frac{k(x-y)}{2}\right)^\ast(t)\right)\leq\left(\frac{k_2}{k_1+k_2}\varphi(k_1x)+\frac{k_1}{k_1+k_2}\varphi(k_2y)\right)^\ast(t).
\end{align*}

In the following we shall show that for some $t\in[0,\alpha)$ there holds
\begin{align*}
\varphi\left(\left(\frac{k(x+y)}{2}\right)^\ast(t)\right)<\left(\frac{k_2}{k_1+k_2}\varphi(k_1x)+\frac{k_1}{k_1+k_2}\varphi(k_2y)\right)^\ast(t)
\end{align*}
or
\begin{align*}
\varphi\left(\left(\frac{k(x-y)}{2}\right)^\ast(t)\right)<\left(\frac{k_2}{k_1+k_2}\varphi(k_1x)+\frac{k_1}{k_1+k_2}\varphi(k_2y)\right)^\ast(t).
\end{align*}

Suppose
\begin{align}\label{3}
\varphi\left(\left(\frac{k(x\pm y)}{2}\right)^\ast(t)\right)&=\left(\varphi\left(\frac{k(x\pm y)}{2}\right)\right)^\ast(t)\notag\\
&=\left(\frac{k_2}{k_1+k_2}\varphi(k_1x)+\frac{k_1}{k_1+k_2}\varphi(k_2y)\right)^\ast(t)
\end{align}
for any $t\in[0,\alpha)$.

Case 2.2.1 Let
\begin{align*}
\left(\frac{k_2}{k_1+k_2}\varphi(k_1x)+\frac{k_1}{k_1+k_2}\varphi(k_2y)\right)^\ast(0)>\left(\frac{k_2}{k_1+k_2}\varphi(k_1x)+\frac{k_1}{k_1+k_2}\varphi(k_2y)\right)^\ast(t)
\end{align*}
for all $t>\alpha$.

Define
\begin{align*}
t_0&=\sup\Bigg\{s:\left(\frac{k_2}{k_1+k_2}\varphi(k_1x)+\frac{k_1}{k_1+k_2}\varphi(k_2y)\right)^\ast(s)\\
&>\left(\frac{k_2}{k_1+k_2}\varphi(k_1x)+\frac{k_1}{k_1+k_2}\varphi(k_2y)\right)^\ast(t) \quad \textit{for\ each}\quad t>\alpha\Bigg\}.
\end{align*}
Obviously, we have $0<t_0\leq\alpha$ and
\begin{align*}
&\left(\frac{k_2}{k_1+k_2}\varphi(k_1x)+\frac{k_1}{k_1+k_2}\varphi(k_2y)\right)^\ast(t_0)\\
&=\left(\frac{k_2}{k_1+k_2}\varphi(k_1x)+\frac{k_1}{k_1+k_2}\varphi(k_2y)\right)^\ast(\alpha)>0.
\end{align*}
If $t_0=\alpha$, then
\begin{align*}
\left(\frac{k_2}{k_1+k_2}\varphi(k_1x)+\frac{k_1}{k_1+k_2}\varphi(k_2y)\right)^\ast(s)>\left(\frac{k_2}{k_1+k_2}\varphi(k_1x)+\frac{k_1}{k_1+k_2}\varphi(k_2y)\right)^\ast(\alpha)
\end{align*}
for any $s<\alpha$, or
\begin{align*}
\left(\frac{k_2}{k_1+k_2}\varphi(k_1x)+\frac{k_1}{k_1+k_2}\varphi(k_2y)\right)^\ast(\alpha)>\left(\frac{k_2}{k_1+k_2}\varphi(k_1x)+\frac{k_1}{k_1+k_2}\varphi(k_2y)\right)^\ast(t)
\end{align*}
for all $t>\alpha$.

If $t_0<\alpha$, there exists $t>\alpha$ such that
\begin{align*}
\left(\frac{k_2}{k_1+k_2}\varphi(k_1x)+\frac{k_1}{k_1+k_2}\varphi(k_2y)\right)^\ast(s)&>\left(\frac{k_2}{k_1+k_2}\varphi(k_1x)+\frac{k_1}{k_1+k_2}\varphi(k_2y)\right)^\ast(t_0)\\
&=\left(\frac{k_2}{k_1+k_2}\varphi(k_1x)+\frac{k_1}{k_1+k_2}\varphi(k_2y)\right)^\ast(t)
\end{align*}
for any $s<t_0$. By (\cite{Krein},Property 7, p.64), we can find the set $e_{t_0}$ with $\mu(e_{t_0})=t_0$ such that
\begin{align*}
\int_{0}^{t_0}\left(\frac{k_2}{k_1+k_2}\varphi(k_1x)+\frac{k_1}{k_1+k_2}\varphi(k_2y)\right)^\ast(t)\,dt\\
=\int_{e_{t_0}}\left(\frac{k_2}{k_1+k_2}\varphi(k_1x(t))+\frac{k_1}{k_1+k_2}\varphi(k_2y(t))\right)\,dt.
\end{align*}
By the proof of (\cite{Krein}, Property 7), we may infer that
\begin{align*}
\frac{k_2}{k_1+k_2}\varphi(k_1x(s))+\frac{k_1}{k_1+k_2}\varphi(k_2y(s))\geq\lim_{t\rightarrow t_0^-}\left(\frac{k_2}{k_1+k_2}\varphi(k_1x)+\frac{k_1}{k_1+k_2}\varphi(k_2y)\right)^\ast(t)
\end{align*}
for $\mu$-a.e $s\in e_{t_0}$. According to the definition of $t_0$, we obtain
\begin{align}\label{4}
&\frac{k_2}{k_1+k_2}\varphi(k_1x(s))+\frac{k_1}{k_1+k_2}\varphi(k_2y(s))\notag\\
&>\left(\frac{k_2}{k_1+k_2}\varphi(k_1x)+\frac{k_1}{k_1+k_2}\varphi(k_2y)\right)^\ast(t)\geq0
\end{align}
for $\mu$-a.e $s\in e_{t_0}$ and each $t>t_0$.
Again using the definition of $t_0$, we get that for $\mu$-a.e. $s\in[0,1)\backslash e_{t_0}$, there exists $t(s)>t_0$ such that
\begin{align}\label{5}
&\frac{k_2}{k_1+k_2}\varphi(k_1x(s))+\frac{k_1}{k_1+k_2}\varphi(k_2y(s))\notag\\
&\leq\left(\frac{k_2}{k_1+k_2}\varphi(k_1x)+\frac{k_1}{k_1+k_2}\varphi(k_2y)\right)^\ast(t(s)).
\end{align}
By equality (\ref{3}), let $e_{t_0}(+)$ and $e_{t_0}(-)$ be sets such that $\mu(e_{t_0}(+))=\mu(e_{t_0}(-))=t_0$ and
\begin{align}\label{6}
&\int_{0}^{t_0}\left(\frac{k_2}{k_1+k_2}\varphi(k_1x)+\frac{k_1}{k_1+k_2}\varphi(k_2y)\right)^\ast(t)\,dt\notag\\
&=\int_{e_{t_0}(+)}\varphi\left(\frac{k(x+y)}{2}(t)\right)\,dt=\int_{e_{t_0}(-)}\varphi\left(\frac{k(x-y)}{2}(t)\right)\,dt.
\end{align}
It is similar to the case of the set $e_{t_0}$, for $\mu$-a.e $s\in e_{t_0}(+)$ and for each $t>t_0$, we get
\begin{align*}
\varphi\left(\frac{k(x+y)}{2}(s)\right)>\left(\frac{k_2}{k_1+k_2}\varphi(k_1x)+\frac{k_1}{k_1+k_2}\varphi(k_2y)\right)^\ast(t).
\end{align*}
By inequality (\ref{4}) and (\ref{5}) we see $e_{t_0}(+)\subset e_{t_0}$. Hence by $\mu(e_{t_0})=t_0=\mu(e_{t_0}(+))$, we get $e_{t_0}(+)=e_{t_0}$. Analogously, we have $e_{t_0}=e_{t_0}(-)$. The equalities (\ref{3}) and (\ref{6}) yield
\begin{align*}
\varphi\left(\frac{k(x+y)}{2}(t)\right)=\varphi\left(\frac{k(x-y)}{2}(t)\right)=\frac{k_2}{k_1+k_2}\varphi(k_1x(t))+\frac{k_1}{k_1+k_2}\varphi(k_2y(t))
\end{align*}
for $\mu$-a.e. $t\in e_{t_0}$. By the convexity of $\varphi$, we know $x(t)=y(t)=0$ for $\mu$-a.e. $t\in e_{t_0}$. Which is a contradiction with (\ref{4}).

Case 2.2.2 Let
\begin{align*}
&\left(\frac{k_2}{k_1+k_2}\varphi(k_1x)+\frac{k_1}{k_1+k_2}\varphi(k_2y)\right)^\ast(0)=\left(\frac{k_2}{k_1+k_2}\varphi(k_1x)+\frac{k_1}{k_1+k_2}\varphi(k_2y)\right)^\ast(\alpha)\\
&=\left(\frac{k_2}{k_1+k_2}\varphi(k_1x)+\frac{k_1}{k_1+k_2}\varphi(k_2y)\right)^\ast(t)>0
\end{align*}
for some $t>\alpha$.

Define
\begin{align*}
A=\Bigg\{t\in[0,1): &\frac{k_2}{k_1+k_2}\varphi(k_1x(t))+\frac{k_1}{k_1+k_2}\varphi(k_2y(t))\\
&=\left(\frac{k_2}{k_1+k_2}\varphi(k_1x)+\frac{k_1}{k_1+k_2}\varphi(k_2y)\right)^\ast(0)\Bigg\},\\
A_+=\Bigg\{t\in[0,1): &\varphi\left(\frac{k(x+y)}{2}(t)\right)\\
&=\left(\frac{k_2}{k_1+k_2}\varphi(k_1x)+\frac{k_1}{k_1+k_2}\varphi(k_2y)\right)^\ast(0)\Bigg\},\\
A_-=\Bigg\{t\in[0,1): &\varphi\left(\frac{k(x-y)}{2}(t)\right)\\
&=\left(\frac{k_2}{k_1+k_2}\varphi(k_1x)+\frac{k_1}{k_1+k_2}\varphi(k_2y)\right)^\ast(0)\Bigg\}.
\end{align*}
The equality (\ref{3}) and  the convexity of $\varphi$ follow that $\mu(A)>\alpha$, $A_+\subset A$, $A_-\subset A$ and $\min\{\mu(A_+), \mu(A_-)\}\geq\alpha$. Let $A_0=A_+\cap A_-=\{t\in A: \min\{|x(t)|,|y(t)|\}=0\}$. Since $\alpha>\frac{1}{2}$, it is easy to see that $\mu(A_0)>0$. According to the definitions of $A$, $A_+$ and $A_-$, we get
\begin{align*}
&\varphi\left(\frac{k(x+y)}{2}(t)\chi_{A_0}(t)\right)=\varphi\left(\frac{k(x-y)}{2}(t)\chi_{A_0}(t)\right)\\
&=\frac{k_2}{k_1+k_2}\varphi(k_1x(t)\chi_{A_0}(t))+\frac{k_1}{k_1+k_2}\varphi(k_2y(t)\chi_{A_0}(t)).
\end{align*}
Hence $x(t)\chi_{A_0}(t)=y(t)\chi_{A_0}(t)=0$ by the convexity of $\varphi$. By virtue of  $A_0\subset A$, we get a contradiction.

From Case 2.2.1 and Case 2.2.2, we infer that there exists $t\in[0,\alpha)$ such that the inequality
\begin{align*}
\varphi\left(\left(\frac{k(x+y)}{2}\right)^\ast(t)\right)<\left(\frac{k_2}{k_1+k_2}\varphi(k_1x)+\frac{k_1}{k_1+k_2}\varphi(k_2y)\right)^\ast(t)
\end{align*}
or
\begin{align*}
\varphi\left(\left(\frac{k(x-y)}{2}\right)^\ast(t)\right)<\left(\frac{k_2}{k_1+k_2}\varphi(k_1x)+\frac{k_1}{k_1+k_2}\varphi(k_2y)\right)^\ast(t)
\end{align*}
holds.
If there exists $t\in[0,\alpha)$ such that
\begin{align*}
\varphi\left(\left(\frac{k(x+y)}{2}\right)^\ast(t)\right)<\left(\frac{k_2}{k_1+k_2}\varphi(k_1x)+\frac{k_1}{k_1+k_2}\varphi(k_2y)\right)^\ast(t),
\end{align*}
then by the right continuity of the rearrangement, we get
\begin{align*}
\rho_{\varphi,\omega}\left(\frac{k(x+y)}{2}\right)&=\int_{0}^{\alpha}\varphi\left(\left(\frac{k(x+y)}{2}\right)^\ast(t)\right)\omega(t)\,dt\\
&<\int_{0}^{\alpha}\left(\frac{k_2}{k_1+k_2}\varphi(k_1x)+\frac{k_1}{k_1+k_2}\varphi(k_2y)\right)^\ast(t)\omega(t)\,dt\\
&=\left\|\frac{k_2}{k_1+k_2}\varphi(k_1x)+\frac{k_1}{k_1+k_2}\varphi(k_2y)\right\|_{1,\omega}\\
&\leq\frac{k_2}{k_1+k_2}(k_1\|x\|^\circ_{\varphi,\omega}-1)+\frac{k_1}{k_1+k_2}(k_2\|y\|^\circ_{\varphi,\omega}-1)\\
&=k-1.
\end{align*}
Similarly, if there exists $t\in[0,\alpha)$ such that
\begin{align*}
\varphi\left(\left(\frac{k(x-y)}{2}\right)^\ast(t)\right)<\left(\frac{k_2}{k_1+k_2}\varphi(k_1x)+\frac{k_1}{k_1+k_2}\varphi(k_2y)\right)^\ast(t),
\end{align*}
we obtain
\begin{align*}
\rho_{\varphi,\omega}\left(\frac{k(x-y)}{2}\right)<k-1.
\end{align*}
In summary, we get
\begin{align*}
\min\left\{\left\|\frac{x+y}{2}\right\|^\circ_{\varphi,\omega}, \left\|\frac{x-y}{2}\right\|^\circ_{\varphi,\omega}\right\}<1.
\end{align*}

\section{Local uniform non-squareness of $\Lambda^\circ_{\varphi,\omega}$}
\begin{theorem}\label{5-th}
Orlicz-Lorentz space $\Lambda^\circ_{\varphi,\omega}[0,\infty)$ is locally uniformly non-square if and only if

(a)\ $\psi\in\Delta_2(\mathbb{R})$,

(b)\ $\int_0^\infty\omega(t)\,dt=\infty$.
\end{theorem}
\proof(Necessity) By Theorem \ref{3-th} and Lemma \ref{5-le}, we can see $\psi\in\Delta_2(\mathbb{R})$ and $\int_0^\infty\omega(t)\,dt=\infty$.

(Sufficiency) By Lemma \ref{1-le}, we see that there exist $\xi_1,\xi_2\in(1,+\infty)$ with $\xi_1<\xi_2$ such that
$k_u\in(\xi_1,\xi_2)$ for any $u\in B(\Lambda^\circ_{\varphi,\omega}[0,\infty))$. Fix $x\in S(\Lambda^\circ_{\varphi,\omega}[0,\infty))$.
For any $y\in B(\Lambda^\circ_{\varphi,\omega}[0,\infty))$. Fix $k_1\in K(x)$, $k_2\in K(y)$. Then
\begin{align*}
\frac{\xi_1}{\xi_1+\xi_2}\leq\frac{k_2}{k_1+k_2}\leq\frac{\xi_2}{\xi_1+\xi_2}.
\end{align*}
Let $k=\frac{2k_1k_2}{k_1+k_2}$. By the range of $k_1$ and $k_2$, we find
\begin{align*}
k=\frac{2k_1k_2}{k_1+k_2}=\frac{2}{\frac{1}{k_1}+\frac{1}{k_2}}\leq\frac{2}{\frac{1}{\xi_2}+\frac{1}{\xi_2}}=\xi_2.
\end{align*}
Certainly there exist $t_1,t_2\in[0,+\infty)$  such that $0\leq t_1<t_2<\infty$,
\begin{align}\label{7}
\int_{t_1}^{t_2}\varphi(k_1x^\ast(t))\omega(t)\,dt\geq\int_{t_1}^{t_2}\varphi(x^\ast(t))\omega(t)\,dt:=\xi\in(0,1]
\end{align}
and $x^\ast(s)>x^\ast(t)>x^\ast(w)$ for all $s\in(0,t_1)$, $t\in(t_1,t_2)$ and $w>t_2$ whenever $t_1>0$,
as well as $x^\ast(t)>x^\ast(w)$ for all $t\in(t_1,t_2)$ and $w>t_2$ whenever $t_1=0$.

By (\cite{Krein}, property 7, page 64), there exist sets $e_{t_1}$ and $e_{t_2}$ with $\mu(e_{t_1})=t_1$ and $\mu(e_{t_2})=t_2$ such that
\begin{align*}
\int_{0}^{t_1}x^\ast(t)\,dt=\int_{e_{t_1}}|x(t)|\,dt \quad \textit{and} \quad \int_{0}^{t_2}x^\ast(t)\,dt=\int_{e_{t_2}}|x(t)|\,dt
\end{align*}
(in the case $t_1=0$ we have $e_{t_1}=\emptyset$). It is easy to see that $e_{t_1}\subsetneq e_{t_2}$. For $u\in[\xi_1x^\ast(t_2),  \xi_2x^\ast(t_1)]$, we can find $\eta\in(0,1)$ such that
\begin{align*}
\varphi\left(\frac{k_2}{k_1+k_2}u\right)\leq (1-\eta)\frac{k_2}{k_1+k_2}\varphi(u).
\end{align*}
Define
\begin{align*}
B_1=\{t\in e_{t_2}\backslash e_{t_1}: x(t)y(t)\geq0\},\\
B_2=\{t\in e_{t_2}\backslash e_{t_1}: x(t)y(t)<0\}.
\end{align*}
By (\ref{7}), we obtain $\rho_{\varphi,\omega}(k_1x\chi_{e_{t_2}\backslash e_{t_1}})\geq\xi$. Obviously, $\rho_{\varphi,\omega}(k_1x\chi_{B_1})\geq\frac{\xi}{2}$ or $\rho_{\varphi,\omega}(k_1x\chi_{B_2})\geq\frac{\xi}{2}$. Assume $\rho_{\varphi,\omega}(k_1x\chi_{B_1})\geq\frac{\xi}{2}$, then we have
\begin{align*}
\varphi\left(\frac{k(x(t)-y(t))}{2}\right)&=\varphi\left(\frac{k_2}{k_1+k_2}k_1x(t)-\frac{k_1}{k_1+k_2}k_2y(t)\right)\\
&\leq\varphi\left(\max\left\{\frac{k_2}{k_1+k_2}k_1|x(t)|,\frac{k_1}{k_1+k_2}k_2|y(t)|\right\}\right)\\
&\leq\varphi\left(\frac{k_2}{k_1+k_2}k_1x(t)\right)+\varphi\left(\frac{k_1}{k_1+k_2}k_2y(t)\right)\\
&\leq(1-\eta)\frac{k_2}{k_1+k_2}\varphi(k_1x(t))+\frac{k_1}{k_1+k_2}\varphi(k_2y(t))
\end{align*}
for $\mu$-a.e. $t\in B_1$. Consequently,
\begin{align*}
\varphi\left(\frac{k(x(t)-y(t))}{2}\right)\leq\frac{k_2}{k_1+k_2}\varphi(k_1x(t))+\frac{k_1}{k_1+k_2}\varphi(k_2y(t))-\frac{\eta k_2}{k_1+k_2}\varphi(k_1x(t)\chi_{B_1}(t)).
\end{align*}
Since $\rho_{\varphi,\omega}(k_1x\chi_{B_1})\geq\frac{\xi}{2}$, we get $\|\eta\varphi(k_1x\chi_{B_1})\|_{1,\omega}\geq\frac{\eta\xi}{2}$. Then by the lower local uniform monotonicity of the Lorentz space $L_{1,\omega}$, we obtain
\begin{align*}
\left\|\varphi\left(\frac{k(x-y)}{2}\right)\right\|_{1,\omega}&\leq\left\|\frac{k_2}{k_1+k_2}\varphi(k_1x)+\frac{k_1}{k_1+k_2}\varphi(k_2y)-\frac{\eta k_2}{k_1+k_2}\varphi(k_1x\chi_{B_1})\right\|_{1,\omega}\\
&\leq\left\|\frac{k_1}{k_1+k_2}\varphi(k_2y)\right\|_{1,\omega}+\left\|\frac{k_2}{k_1+k_2}\varphi(k_1x)-\frac{\eta k_2}{k_1+k_2}\varphi(k_1x\chi_{B_1})\right\|_{1,\omega}\\
&\leq\frac{k_1}{k_1+k_2}\|\varphi(k_2y)\|_{1,\omega}+\frac{k_2}{k_1+k_2}\|\varphi(k_1x)\|_{1,\omega}-\frac{k_2}{k_1+k_2}\delta_1,
\end{align*}
where $\delta_1=\delta_1\left(\varphi(k_1x),\frac{\eta\xi}{2}\right)$ is a constant only depends on $x$.
That is to say
\begin{align*}
\rho_{\varphi,\omega}\left(\frac{k(x-y)}{2}\right)&\leq\frac{k_1}{k_1+k_2}\rho_{\varphi,\omega}(k_2y)+\frac{k_2}{k_1+k_2}\rho_{\varphi,\omega}(k_1x)-\frac{k_2}{k_1+k_2}\delta_1\\
&\leq\frac{k_1}{k_1+k_2}(k_2-1)+\frac{k_2}{k_1+k_2}(k_1-1)-\frac{k_2}{k_1+k_2}\delta_1\\
&\leq k-1-\frac{k_2}{k_1+k_2}\delta_1.
\end{align*}
By the definition of Orlicz norm, we get
\begin{align*}
\left\|\frac{x-y}{2}\right\|^\circ_{\varphi,\omega}&\leq\frac{1+\rho_{\varphi,\omega}(\frac{k(x-y)}{2})}{k}\\
&\leq1-\frac{k_2}{k(k_1+k_2)}\delta_1\\
&\leq1-\delta,
\end{align*}
where $\delta=\frac{\xi_1}{\xi_2(\xi_1+\xi_2)}\delta_1$ only depends on $x$.

If $\rho_{\varphi,\omega}(k_1x\chi_{B_2})\geq\frac{\xi}{2}$, we can prove $\left\|\frac{x+y}{2}\right\|^\circ_{\varphi,\omega}\leq1-\delta$ analogously.

\begin{theorem}\label{6-th}
Let $\alpha=\sup\{t\geq0: \omega(t)>0\}=1$, Then Orlicz-Lorentz space $\Lambda^\circ_{\varphi,\omega}[0,1)$ is locally uniformly non-square if and only if $\psi\in\Delta_2(\mathbb{\infty})$.
\end{theorem}
\proof(Necessity) By Lemma \ref{5-le}, it has been proved.

(Sufficiency) By Lemma \ref{1-le}, we see that there exist $\xi_1,\xi_2\in(1,+\infty)$ with $\xi_1<\xi_2$ such that
$k_u\in(\xi_1,\xi_2)$ for any $u\in B(\Lambda^\circ_{\varphi,\omega}[0,1))$. Fix $x\in S(\Lambda^\circ_{\varphi,\omega}[0,1))$.
For any $y\in B(\Lambda^\circ_{\varphi,\omega}[0,1))$. Fix $k_1\in K(x)$, $k_2\in K(y)$. Let $k=\frac{2k_1k_2}{k_1+k_2}$ and $A_x=\{t\in[0,1):|x(t)|>0\}$. Obviously, $\mu(A_x)>0$. Hence, there exist $t_1, t_2$ such that $0\leq t_1<t_2\leq t_0$, $0<x^\ast(t_2)\leq x^\ast(t_1)<\infty$ and
\begin{align*}
\int_{t_1}^{t_2}\varphi(k_1x^\ast(t))\omega(t)\,dt\geq\int_{t_1}^{t_2}\varphi(x^\ast(t))\omega(t)\,dt>0.
\end{align*}
Analogously as in the proof of sufficiency in Theorem \ref{5-th}, we can find $\delta>0$ only depends on $x$ such that
\begin{align*}
\min\left\{\left\|\frac{x-y}{2}\right\|^\circ_{\varphi,\omega},\left\|\frac{x+y}{2}\right\|^\circ_{\varphi,\omega}\right\}\leq1-\delta.
\end{align*}

\end{document}